\documentclass{article}
\usepackage[russian]{babel}
\usepackage{mathrsfs}
\usepackage[usenames]{color}
\usepackage[dvips]{graphicx}

\begin{document}
\pagestyle{myheadings}
\title{High performance parallel algorithm for solving elliptic equations
with non-separable variables.}
\author{\sc{А.В. Терехов} \\ \\
{Institute of Computational Mathematics and Mathematical
Geophysics,
630090,Novosibirsk,Russia}\\
{Novosibirsk State University, 630090, Novosibirsk, Russia} \small
\it} \maketitle
\begin{abstract}

A parallel algorithm for computing the finite difference solution
to the elliptic equations with non-separable variables is
presented. The resultant matrix is symmetric positive definite,
thus the preconditioning conjugate gradient or the chebyshev
method can be applied. A differential analog to the Laplace
operator is used as preconditioner. For inversion of the Laplace
operator we implement a parallel version of the separation
variable method, which includes the sequential FFT algorithm and
the parallel solver for tridiagonal matrix equations (dichotomy
algorithm). On an example of solving acoustic equations by the
integral Laguerre transformation  method, we show that the
algorithm proposed is highly efficient for a large number of
processors.
\\\\
\textbf{Ключевые слова}: Трехдиагональные системы уравнений,
параллельный алгоритм прогонки, Parallel Dichotomy Algorithm,
метод разделения переменных, эллиптические уравнения с
неразделяемыми переменными, преобразование Лагерра, уравнение
акустики.
\end{abstract}

\section{Введение.}
Существует обширный класс задач (моделирование волновых полей,
полупроводниковых приборов, задачи теплообмена, гидродинамики и
др.) \cite{Samarsk_Teplo1,Andreson,Pantakar,Hockney:Istwood}, в
рамках которых требуется решить не одну, а серию уравнений с
постоянным оператором и различными правыми частями

\begin{equation}
 \label{main_problem}\nabla \left[\kappa({\bf r})\, \nabla
u_i\right]-q({\bf r})u_i=-f_i({\bf r}),\quad{\bf r} \in
\mathbf{R}^n,\; i=1,2,3,\quad \kappa({\bf
r})>0,\;q(\mathbf{r})\geq0\quad i=1,...M .
\end{equation} В силу того, что при решении современных проблем
математического моделирования количество задач в серии может
достигать нескольких тысяч, разумно проводить такие расчеты на
суперкомпьютере.

К сожалению, для большинства экономичных численных методов,
разработанных для решения разностных уравнений
(попеременно-треугольный метод, метод переменных направлений,
циклической редукции  и т.д.) \cite{Samarski_Nikolaev}, их
реализация на многопроцессорных вычислительных системах не
приводит к значительному сокращению времени счета. Это является
следствием того, что при теоретических оценках качества
последовательных алгоритмов ограничиваются только подсчетом числа
арифметических операций, необходимых для отыскания решения с
заданной точностью. В то же время анализ работ по данной тематике
\cite{Dongarra:Linear,Dongarra:Reduce} показывает, что на этапе
конструирования численного метода для решения задач на
суперкомпьютере необходимо дополнительно учитывать временные
затраты на межпроцессорные взаимодействия. Несмотря на то что
решению задачи (\ref{main_problem}) на суперкомпьютере посвящено
значительное число работ
\cite{Ell_prob_4,Ell_prob_5,Ell_prob_1,Ell_prob_2,Ell_prob_3,Ell_prob_6},
есть обстоятельства, которые делают исследования в этой области
весьма актуальными. Так, постоянный рост числа (процессорных
элементов) ПЭ, объединенных в рамках одного суперкомпьютера,
предъявляет все более жесткие требования к масштабируемости
численных алгоритмов. Например, алгоритм циклической редукции
\cite{Samarski_Nikolaev,Ilin_Kuznecov,hockney:cyclic,Johnsson:cyclic}
эффективен при использовании относительно небольшого числа
процессоров, поэтому его реализация на современных параллельных
вычислительных системах сильно ограничена. Таким образом, одним из
требований, предъявляемых к разрабатываемому параллельному
алгоритму, является возможность его эффективной реализации на
суперкомпьютере для большого числа процессоров (от тысячи и
более).

Дополнительно от метода решения задачи (\ref{main_problem})
необходимо требовать его экономичности, так как в противном случае
смысл параллельной реализации может быть в значительной степени
утрачен.

В продолжении работы \cite{Terekhov:dichotomy}, в которой
предложен новый алгоритм параллельной прогонки для
мультикомпьютера: Алгоритм Дихотомии, в статье исследован
параллельный алгоритм для решения серии уравнений вида
(\ref{main_problem}) в цилиндрической системе координат на
прямоугольной сетке. Эффективность предлагаемого подхода
исследована в рамках реализации численно-аналитического метода
решения уравнения акустики \cite{Mikh:Laguerre}.

Структура статьи следующая. В части 2 формулируется краевая задача
для эллиптического уравнения второго порядка в
аксиально-симметричной геометрии. Обсуждается вопрос о выборе
итерационной процедуры решения разностных уравнений с учетом ее
последующей реализации на суперкомпьютере. В качестве
предобусловливающего оператора рассматривается оператор Лапласа. В
части 3 рассмотрен алгоритм параллельной прогонки (Алгоритм
Дихотомии). В части 4 рассматривается модельная задача о
распространении акустических волн в неоднородной среде. Приводятся
оценки точности получаемого решения, исследуется эффективность
распараллеливания. В части 5 подводятся итоги проделанной работы.
\section{Параллельный Алгоритм решения разностных уравнений.}
\subsection{Постановка задачи.}
В прямоугольной области $G=\{0\leq r \leq l_1,\; 0\leq z \leq
l_2\}$ с границей $\Gamma$ рассмотрим следующую краевую задачу
вида (\ref{main_problem})

\begin{equation}
\left\{
\begin{array}{ll}
\label{main_eq} \displaystyle \frac{1}{r}\frac{\partial}{\partial
r}\left(r
\kappa(r,z)\frac{\partial u}{\partial r}\right)+\frac{\partial}{\partial z}\left(\kappa(r,z)\frac{\partial u}{\partial z}\right) -q(r,z)u=-f(r,z),& (r,z) \in G,\\\\
\displaystyle \left. r\kappa\frac{\partial u}{\partial
r}\right|_{r=0}=\left.\frac{\partial u}{\partial
z}\right|_{z=0,l_2}=\left.u\right|_{r=l_1}=0,&(r,z) \in \Gamma.
\end{array}\right.
\end{equation}

Для обеспечения существования и единственности решения задачи
(\ref{main_eq}) потребуем выполнения следующих неравенств
\cite{Samarski_Nikolaev}:
\begin{equation}
\label{elliptic_condition}
\begin{array}{ll}
0<s_1 \leq \kappa(r,z) \leq s_2,\\
0<d_1 \leq q(r,z)\leq d_2.
\end{array}
\end{equation}

На прямоугольной сетке $ \label{grid2}
\bar{\omega}=\bar{\omega}_r\times\bar{\omega}_z=\omega\bigcup\gamma$
,где
$$
\begin{array}{l} \bar \omega_r=\left\{r_i=(i-0.5)h_1,\;
i=1,...,N_1,\; h_1=l_1/(N_1-0.5)\right\},\\\\
\bar \omega_z=\left\{z_k=(k-0.5)h_2,\; k=1,...,N_2,\;
h_2=l_2/(N_2-0.5) \right\}
\end{array}
$$
$\omega=\bar\omega\bigcap G,\quad \gamma=\bar\omega\bigcap\Gamma$
дифференциальной задаче (\ref{main_eq}) соответствует разностная
задача\cite{Samarski_Nikolaev,Samarski_andreev}

\begin{equation}
\begin{array}{ll}
\displaystyle\label{poisson_diff_2}
\left(\Lambda_r+\Lambda_z\right)y-w(x)y=-\phi(x),\quad x \in
\bar{\omega},
\end{array}
\end{equation}

\begin{equation}
\label{scheme_op_r_z}
\begin{array}{lr}
\Lambda_ry=\left\{
\begin{array}{ll}
\displaystyle \frac{1}{h_r}a_1y_{r},& i=1\\\\
\displaystyle \left(a_1y_{\bar{r}}\right)_{r},&1\leq i\leq N_1-1
\end{array}\right., &
\Lambda_zy=\left\{
\begin{array}{ll}
\displaystyle \frac{1}{h_z}a_2y_{z},& k=1\\\\
\displaystyle \left(a_2y_{\bar{z}}\right)_{z},&1\leq k\leq
N_2-1\\\\
\displaystyle \frac{1}{h_z}a_2y_{\bar{z}},&k=N_2
\end{array}.\right.
\end{array}
\end{equation}

Коэффициенты разностной схемы будем определять по формулам
\begin{equation}
\label{diff_approx1}
\begin{array}{ll}
 \displaystyle
a_1(i,k)=\bar{r}_i\kappa\left(\bar{r}_i,z_k\right),&
a_2(i,k)=r_i\kappa\left(r_i,\bar{z}_k\right)\\
w(i,k)=r_iq(r_i,z_k),& \phi(i,j)=r_if(r_i,z_k), %\\
\end{array}
\end{equation}
 где $\bar{r}_i=r_i+0.5h_1,\; \bar{z}_k=z_k+0.5h_2$; $\;y_{\bar r},\,y_{\bar
 z}$ и $y_r,\,y_z$ -- разностные соотношения по
 $z$ и по $r$ "назад"\ и "вперед"\
 \cite{Samarski_andreev}.
Краевое условие на стороне $r=l_2$ аппроксимируем точно $
y_{N_1,k}=0, \quad k=1,...,N_2$.

Схема (\ref{poisson_diff_2})-(\ref{diff_approx1}) на сетке
$\bar{\omega}$ обладает вторым порядком аппроксимации и является
консервативной. Поскольку для схемы (\ref{poisson_diff_2})
проблемы точности, сходимости и устойчивости достаточно полно
изучены например в \cite{Samarski_andreev}, в дальнейшем мы не
будем останавливаться на этих вопросах, а обратим основное
внимание на параллельные алгоритмы решения разностного уравнения
(\ref{poisson_diff_2}).

\subsection{Алгоритм решения разностного уравнения.}
Разностную задачу (\ref{poisson_diff_2}) будем трактовать как
операторное уравнение в вещественном конечномерном гильбертовом
пространстве $H$
\begin{equation}
\label{a-eq}
 Au=f,\quad A:H \rightarrow H,
\end{equation}
где $A$ не вырожденный, линейный, самосопряженный
$\left(A=A^{*}\right)$ и положительно определенный оператор
$\left(A>0 \right)$ \cite{Samarski_Nikolaev}.

Для приближенного решения задачи (\ref{a-eq}) рассмотрим неявный
итерационный процесс

\begin{equation}
\begin{array}{cl}
\label{main-scheme} \displaystyle
B\frac{\bar{y}-y_{k}}{\tau_{k+1}}+Ay_k=f, \quad B:H \rightarrow H,\\\\
\displaystyle
y_{k+1}=\alpha_{k+1}\bar{y}+(1-\alpha_{k+1})y_{k-1},& k=1,2...
\end{array}
\end{equation}
с произвольным начальным приближением $y_0 \in H$
\cite{Samarski_Nikolaev}. Параметры $\{\tau_k\}$ и оператор $B$
(предобуславливающий) следует выбирать из условия минимума числа
итераций, при котором

\begin{equation}
\label{ocenka} \parallel z_n\parallel = \parallel y_n-u
\parallel_{} \leq \varepsilon\parallel y_0-u
\parallel_{},
\end{equation}
где $u$ -- точное решение задачи (\ref{a-eq}), а $\varepsilon>0$
-- требуемая точность.

Рассмотрим вопрос о выборе предобуславливающего оператора $B$ из
класса операторов, обладающими следующими свойствами:
\begin{enumerate}
    \item $B=B^*>0$, самосопряженность и положительная определенность.
    \item Оператор $B$ должен быть "энергетический эквивалентен"\,
    оператору $A$ в смысле неравенств
    \begin{equation}
    \label{energy_rel}
    \gamma_1 \left(Bu,u\right) \leq \left(Au,u\right) \leq
    \gamma_2 \left(Bu,u\right), \forall u;  \quad \quad
    0<\gamma_1\leq\gamma_2,
    \end{equation}
     где $$\gamma_1=\min_{x\neq
     0}\frac{\left(Ax,x\right)}{\left(Bx,x\right)},\quad
     \gamma_2=\max_{x\neq 0}\frac{\left(Ax,x\right)}{\left(Bx,x\right)}.$$
    \item Оператор $B$ должен быть легко обратимым, по сравнению с
    оператором $A$.
    \item Так как рассматриваемый алгоритм предполагается
реализовывать на многопроцессорных вычислительных системах,
естественно требовать, чтобы процедура обращения
предобуславливающего оператора допускала эффективную параллельную
реализацию. \label{x1}
\end{enumerate}

В работах
\cite{Samarski_Nikolaev,bernhardt,Concus:Golub,Israeli,Chang,Waveform}
исследованы различные алгоритмы для решения эллиптических
уравнений второго порядка с переменными коэффициентами, в которых
итерационный процесс сводится к многократному решению уравнения
Пуассона. Известно\cite{Samarski_Nikolaev}, что случае
$B\equiv\Delta$ оценка числа обусловленности
$cond\left(B^{-1}A\right)=\gamma_2/\gamma_1$ не зависит от шага
сетки, поэтому данный класс предобуславливателей эффективен,
например, в задачах сейсмической разведки \cite{White}, где при
умеренной контрастности среды требуется проводить расчеты с
большим числом узлов, приходящимся на характерную длину волны.

В \cite{Terekhov:dichotomy} предложен алгоритм дихотомии для
решения трехдиагональный систем уравнений, который позволяет
достичь высокой производительности при параллельной реализации
таких алгоритмов, как метод разделения переменных
\cite{Hockney:Istwood,Samarski_Nikolaev,Hockney}, метод переменных
направлений \cite{Samarski_Nikolaev,Ianenko,Pissman}, циклической
редукции \cite{Hockney:Istwood,Samarski_Nikolaev,Bunemann}.
Использование этих параллельных процедур для обращения разностного
аналога оператора Лапласа, позволяет удовлетворить требованию 4,
предъявляемому к предобуславливающему оператору и выбрать его в
виде\cite{Waveform}
\begin{equation}
 \label{coeff2}
B\equiv \frac{s_1+s_2}{2}\Delta-\frac{d_1+d_2}{2}.
\end{equation}

Поскольку $AB^{-1}A=\left(AB^{-1}A\right)^{*}$, то для вычисления
итерационных параметров можно использовать метод Чебышева или
Метод сопряженных градиентов\cite{Samarski_Nikolaev}. С точки
зрения параллельной реализации преимущество метода Чебышева по
сравнению с вариационными методами состоит в отсутствие операции
скалярного произведения над распределенными данными, а
следовательно, и коммуникационных взаимодействий типа
"All-to-All-Reduce"\ \cite{Voevodin_MPI,MPI_REDUCTION}. Наличие
такого рода коллективных взаимодействий снижает производительность
и масштабируемость параллельного алгоритма\cite{Dongarra:Linear}.

    В случае высококонтрастных сред выбранная предобуславливающая
процедура может и не обеспечивать высокой скорости сходимости
итерационного процесса, более того, иногда наблюдается его
расходимость\cite{Waveform}. В этом случае можно рекомендовать
методы с выделением границы разрывов сред, которые также требуют
обращения оператора Лапласа\cite{LI1,LI2}.

\section{Алгоритм Дихотомии для решения серии трехдиагональных систем уравнений.}
 Использование таких методов, как метод разделения переменных,
метод переменных направлений или циклической редукции в рамках
схемы (\ref{main-scheme}) для обращения оператора $B\equiv\Delta$,
требует многократного решения трехдиагональных СЛАУ с одной и той
же матрицей, но различными правыми частями

\begin{equation}
\label{main_eq3}
 A {\bf X_{i}}= { \bf F_{i}} ,\quad i=1,2...M,
\end{equation}

$$
A=\left\|%
\begin{array}{ccccc}
  b_1 & a_1 &  &    &  \Large 0\\
  c_2 & b_2 & a_2   &  &  \\
   %& a_2 & b_3 & a_3 &  &  \\
    & \ddots & \ddots & \ddots &  \\
     &  & c_{n-1} & b_{n-1} & a_{n-1} \\
   0  &  &  & c_{n} & b_n \\
\end{array}%
\right\|.
$$

Рассмотрим Алгоритм Дихотомии \cite{Terekhov:dichotomy},
разработанный для решения этой проблемы.
\subsection{Декомпозиция данных.}
Выбор способа декомпозиции данных задачи оказывает существенное
влияние на вычислительную и коммуникационную трудоемкость
параллельного алгоритма, а в итоге и на время счета\cite{Malysh}.
Поэтому рассмотрим вопрос отображения данных задачи
(\ref{main_eq3}) на множество процессорных элементов (ПЭ).

Пусть число ПЭ равно $p$. Разобьем вектор правой части и вектор
решения ${ \bf F}$ и ${ \bf X}$ на подвекторы  ${\bf Q_{i},\,
U_{i}}$ следующим образом:

\begin{equation}
{\bf F}=\left({ \bf Q}_{1}, { \bf Q}_{2},...,{\bf Q}_{p}
\right)^\mathrm{T}=\left(f_1,f_2,...,f_{\it{size}\{\bf{
F}\}-1},f_{\it{size}\{{\bf F}\}}\right)^\mathrm{T
}\label{decom_f},
\end{equation}

\begin{equation}
{\bf X}=\left({\bf U}_{1},{\bf U}_{2},...,{\bf U}_{p}
\right)^\mathrm{T}=\left(x_1,x_2,...,x_{\it{size}\{\bf{X}\}-1},x_{\it{size}\{{\bf
X}\}}\right)^\mathrm{T}. \label{decom_x}
\end{equation}

Полагая число элементов некоторого вектора $\bf V$ равным
$\it{size}\{\bf V\}$, потребуем выполнения следующих условий

$$
\begin{array}{l}
 size\{{\bf Q}_{i}\}=size\{{\bf U}_{i}\} \geq 2 \quad i=1,...,p,\\\\
\sum_{i=1}^{p}{\it size}\{{\bf Q}_{i}\}=\sum_{i=1}^{p}{\it
size}\{{\bf U}_{i}\}={\it size}\{{\bf X}\}.
\end{array}
$$

Будем считать, что ПЭ с номером $i$ принадлежит пара подвекторов
$\left({\bf Q_{i}},{\bf U_{i}}\right)$, а строка матрицы $A$ с
номером $j$ принадлежит тому ПЭ, на котором расположена пара
элементов $\left(x_j,f_j\right)$ из
(\ref{decom_f}),(\ref{decom_x}).

Дополнительно введем следующие обозначения:
\begin{itemize}
    \item Обозначим первый и последние компоненты
некоторого вектора $\bf V$  как $first\{{\bf V}\}$ и $last\{{\bf
V}\}$.

    \item Примем за $\left\{A\right\}_{l}^{t}$ матрицу, получающуюся из
матрицы $A$ путем отбрасыванием всех строк и столбцов с номерами
меньшими $l$ или большими $t$.
    \item Примем за $\left\{\mathbf{V}\right\}_{l}^{t}$ подвектор, получающийся из
вектора $\mathbf{V}$ путем отбрасыванием элементов с номерами
меньшими $l$ или большими $t$.

\end{itemize}

\subsection{Дихотомия СЛАУ.}
\begin{figure}[!htb]  \center
\includegraphics[width=0.5\textwidth]{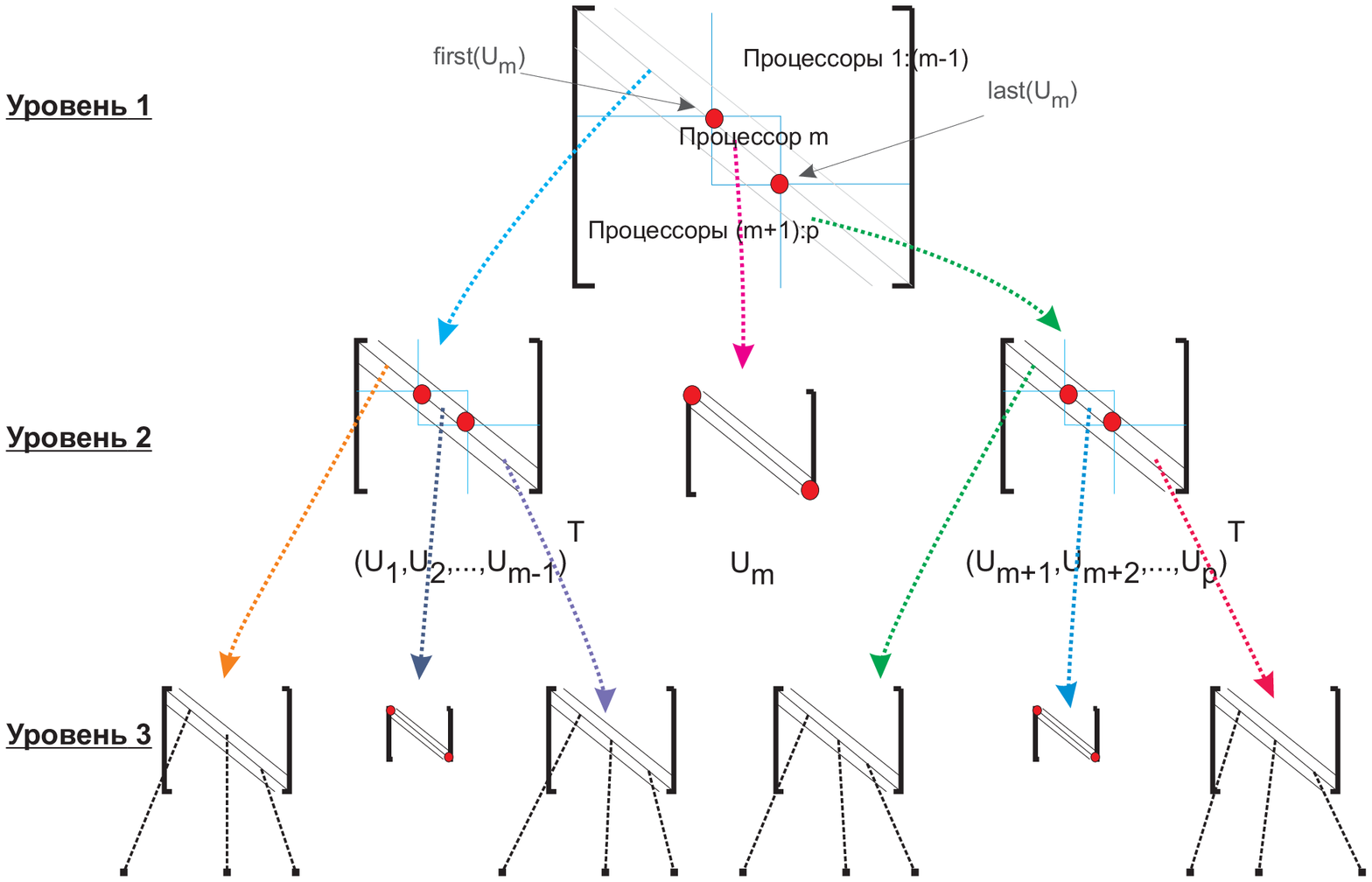}
\caption{Процесс разделения трехдиагональной СЛАУ на независимые
подсистемы. } \label{pic:example_dichotomy}
\end{figure}

Алгоритм Дихотомии \cite{Terekhov:dichotomy} является
представителем класса алгоритмов известного, как "Divide \&
Conquer" \cite{Konovalov,Kon2,Wang}. На каждом уровне дихотомии
(рис.~\ref{pic:example_dichotomy}) трехдиагональная система
уравнений, полученная на предыдущем шаге, путем вычисления решений
в \mbox{$\left(\mathbf{X}\right)_{m_L}=first\{{\bf U}_{m}\}$},
\mbox{$\left(\mathbf{X}\right)_{m_R}=last\{{\bf U}_{m}\}$} --
компонентах разделяется на три независимых подсистемы меньших
размерностей

\begin{equation}\label{systema}
\left\{A\mathbf{X}\right\}_{1}^{m_L-1}=\left\{\mathbf{F}\right\}_{1}^{m_L-1}-a_{m_L-1}first\{\mathbf{U}_m\}\mathbf{e}^{\mathrm{L}},
\end{equation}

\begin{equation}\label{systemb}
\left\{A\mathbf{X}\right\}_{m_L+1}^{m_R-1}=\left\{\mathbf{F}\right\}_{m_L+1}^{m_R-1}-c_{m_L+1}\mathbf{e}^{\mathrm{R}}first\{\mathbf{U}_m\}-a_{m_R-1}\mathbf{e}^{\mathrm{L}}last\{\mathbf{U}_m\},
\end{equation}

\begin{equation}\label{systemc}
\left\{A\mathbf{X}\right\}_{m_R+1}^{n}=\left\{\mathbf{F}\right\}_{m_R+1}^{n}-c_{m_R+1}last\{\mathbf{U}_m\}\mathbf{e}^{\mathrm{R}},
\end{equation}

$$
\mathbf{e}^\mathrm{R}=\left(1,0,0,...,0\right)^{\mathrm{T}},\;\mathbf{e}^\mathrm{L}=\left(0,...,0,0,1\right)^{\mathrm{T}}.
$$

Далее, алгоритм разделения рекурсивно применяется к подсистемам
(\ref{systema}),(\ref{systemc}). В итоге через $\lceil \log_2p
\rceil$ шагов исходная система уравнений (\ref{main_eq3}) будет
разделена на $p$ независимых подсистем вида (\ref{systemb}). На
заключительном шаге Алгоритма Дихотомии на каждом процессоре
определяется решение задач (\ref{systemb}) с помощью любого
последовательного варианта алгоритма прогонки
\cite{Samarski_Nikolaev,Ilin_Kuznecov}.

Достаточным условием применимости Алгоритм Дихотомии является
наличие диагонального преобладания матрицы СЛАУ. По точности,
числу арифметических операций и количеству коммуникационных
взаимодействий Алгоритм Дихотомии  практически эквивалентен методу
циклической редукции \cite{Ilin_Kuznecov}. Однако при сопоставимых
объемах передаваемых данных реальное время межпроцессорных
взаимодействий Алгоритм Дихотомии существенно меньше. Это
объясняется тем, что все межпроцессорные обмены могут быть
осуществлены через последовательность вызовов коллективной
операции "All-to-One-Reduce(+)"\ . Учет таких свойств операции
"$+$"\ , как ассоциативность и коммутативность, позволяет
уменьшить время коммуникационных взаимодействий за счет
возможности их
оптимизации\cite{Tuning:1,Tuning:2,Tuning:3,Collective:1}.
Существенным является тот факт, что оптимизация межпроцессорных
взаимодействий может быть произведена на уровне коммуникационной
библиотеки или языка программирования, что обеспечивает высокую
межплатформенную переносимость программного обеспечения,
включающего в себя Алгоритм Дихотомии.

\subsection{Вычисление $first,last$--компонент.}
Процесс вычисления  $first\{{\bf U}_{i}\}$, $\,last\{{\bf
U}_{i}\}$ компонент состоит из двух этапов: подготовительного,
который выполняется один раз для всех правых частей
(\ref{main_eq3}) и процесса дихотомии, на котором вычисляются
решения для каждой правой части.
\subsubsection{Подготовительные вычисления.}
На подготовительном этапе на $m$-ом процессоре локально без
коммуникационных взаимодействий вычисляются две строки матрицы
$A^{-1}$
\begin{equation}
\begin{array}{cc}
\label{green_v} A^{\mathrm{T}} {\bf G}^{\mathrm{L}}_m={\bf
e}^{}_{m_L}, & A^{\mathrm{T}}{\bf G}^{\mathrm{R}}_m={\bf
e}^{}_{m_R},
\end{array}
\end{equation}

здесь и далее  $m_L\,,m_R$ определены как
\begin{equation}
\label{mrml} first\{{\bf U}_{m}\}=\left(\mathbf{X}\right)_{m_L} ,
last\{{\bf U}_{m}\}=\left(\mathbf{X}\right)_{m_R},
\end{equation} a $\mathbf{e}_k$-- орт
-вектор.

Дополнительно на подготовительном этапе вычисляются два вектора
\begin{equation}
\label{z_vector}
\begin{array}{l}
\mathbf{Z}_{m}^\mathrm{L}=\left(z^{\mathrm{L}}_{1},z^{\mathrm{L}}_{2},...,z^{\mathrm{L}}_{m_L-1},1\right)^{\mathrm{T}},\\\\
\mathbf{Z}_{m}^\mathrm{R}=\left(1,z^{\mathrm{R}}_{m_R+1},\dots,z^{\mathrm{R}}_{n-1},z^{\mathrm{R}}_{n}\right)^{\mathrm{T}},
\end{array}
\end{equation}

компоненты которых определяются из решения систем
\begin{equation} \left\{A\right\}_{1}^{m_L-1}\left(\begin{array}{l}
z^{\mathrm{L}}_1
\\ z^{\mathrm{L}}_2
\\ \dots
%\\ z^{\mathrm{L}}_{m_L-2}
\\ z^{\mathrm{L}}_{m_L-1}
\end{array}\right)=\left(\begin{array}{c}
0
%\\ 0
\\ \dots
\\ 0
\\ -a_{m_L-1}\end{array}\right),\\\\
\label{zl}
\end{equation}

\begin{equation}
\left\{A\right\}_{m_R+1}^{n}\left(\begin{array}{l}
z^{\mathrm{R}}_{m_R+1}
\\ z^{\mathrm{R}}_{m_R+2}
\\ \dots
%\\ z^{\mathrm{R}}_{n-1}
\\ z^{\mathrm{R}}_{n}
\end{array}\right)=\left(\begin{array}{c}-c_{m_R+1}
\\ 0
\\ \dots
%\\ 0
\\ 0
\end{array}\right).\\\\
\label{zr}
\end{equation}

Затраты на подготовительные вычисления для Алгоритма Дихотомии
составят $size\{\mathbf{Z}_m^{\mathrm{R,L}}\}=O(n).$ Таким
образом, Алгоритм Дихотомии имеет смысл применять для решения
нескольких СЛАУ с одной и той же матрицей (\ref{main_eq3}) и
различными правыми частями, т.е.  в случае, когда
подготовительными вычислениями можно пренебречь. Отметим, что
объем вспомогательных вычислений Алгоритма Дихотомии больше по
сравнений с алгоритмом \cite{Konovalov,Kon2,Wang}, тем не менее
эти затраты окупятся за счет сокращения времени коммуникаций.

\subsubsection{Основная формула.}
На втором этапе Алгоритма Дихотомии для вычисления $firts,last$ --
элементов используется представление\cite{Terekhov:dichotomy}
\begin{equation}
\label{theor_3} \left(\mathbf{X}\right)_k=\left\{
\begin{array}{ll}\displaystyle
\sum_{j=1}^{m-1}
\beta_{j}^\mathrm{R}\left(\mathbf{Z}_j^\mathrm{R}\right)_k+\sum_{j=m}^{p}
\beta^\mathrm{L}_{j}\left(\mathbf{Z}_j^\mathrm{L}\right)_k, &
\left(\mathbf{X}\right)_k=first\{\mathbf{U}_m\},\\
\\\\
\displaystyle   \sum_{j=1}^{m}
\beta_{j}^\mathrm{R}\left(\mathbf{Z}_j^\mathrm{R}\right)_k+\sum_{j=m+1}^{p}
\beta^\mathrm{L}_{j}\left(\mathbf{Z}_j^\mathrm{L}\right)_k, &
\left(\mathbf{X}\right)_k=last\{\mathbf{U}_m\},
\end{array}\right.
\end{equation}

где

\begin{equation}
\label{beta_symn}
\begin{array}{ll}
\displaystyle
\beta_m^\mathrm{L}=\sum_{j=m_L}^{m_R}\left(\mathbf{F}\right)_j\left(\mathbf{G}^\mathrm{L}_{m}\right)_j,&
\displaystyle
\beta_m^\mathrm{R}=\sum_{j=m_L}^{m_R}\left(\mathbf{F}\right)_j\left(\mathbf{G}^\mathrm{R}_{m}\right)_j,
\end{array}
\end{equation}
здесь индексы $m_R,m_L$ определяются локально на каждом процессоре
согласно (\ref{mrml}).

Таким образом, при реализации Алгоритма Дихотомии расчет
$first,last$ -- компонент сводится к вычислению сумм
(\ref{theor_3}), тогда как алгоритмы \cite{Konovalov,Kon2,Wang}
требуют решения вспомогательной СЛАУ размерности $2p-2$. Очевидно,
что вычисление сумм на многопроцессорной вычислительной системе
может быть реализовано с большей эффективностью, чем алгоритм
исключения Гаусса, поэтому Алгоритм Дихотомии позволяет достичь
более высокой производительности по сравнению с методами
\cite{Konovalov,Kon2,Wang} для задач вида (\ref{main_eq3}).

\subsubsection{MPI--реализация процесса дихотомии.}
Рассмотрим последовательность действий для решения
трехдиагональной СЛАУ с помощью Алгоритма Дихотомии на $p=7$
процессорах.

В соответствии со спецификацией (\ref{decom_f}), (\ref{decom_x}),
на процессоре с номером $m$ распределяются подвекторы
$\left(\mathbf{U}_m,\mathbf{Q}_m\right)$. На подготовительном
этапе, который выполняется один раз для всех правых частей, на
$m$--ом процессоре для определения векторов
$\mathbf{Z}_m^{\mathrm{R}}$,$\mathbf{Z}_m^{\mathrm{L}}$,$\mathbf{G}_m^{\mathrm{R}}$,$\mathbf{G}_m^{\mathrm{L}}$
решаются задачи (\ref{green_v}),(\ref{zl}),(\ref{zr}) .

В начале второго этапа на процессоре с номером $m$ вычисляются
локальные величины $\beta_m^{\mathrm{R}},\;\beta_m^{\mathrm{L}}$ в
соответствии с (\ref{beta_symn}).

\begin{figure}[!h]
\center
\includegraphics[width=0.8\textwidth]{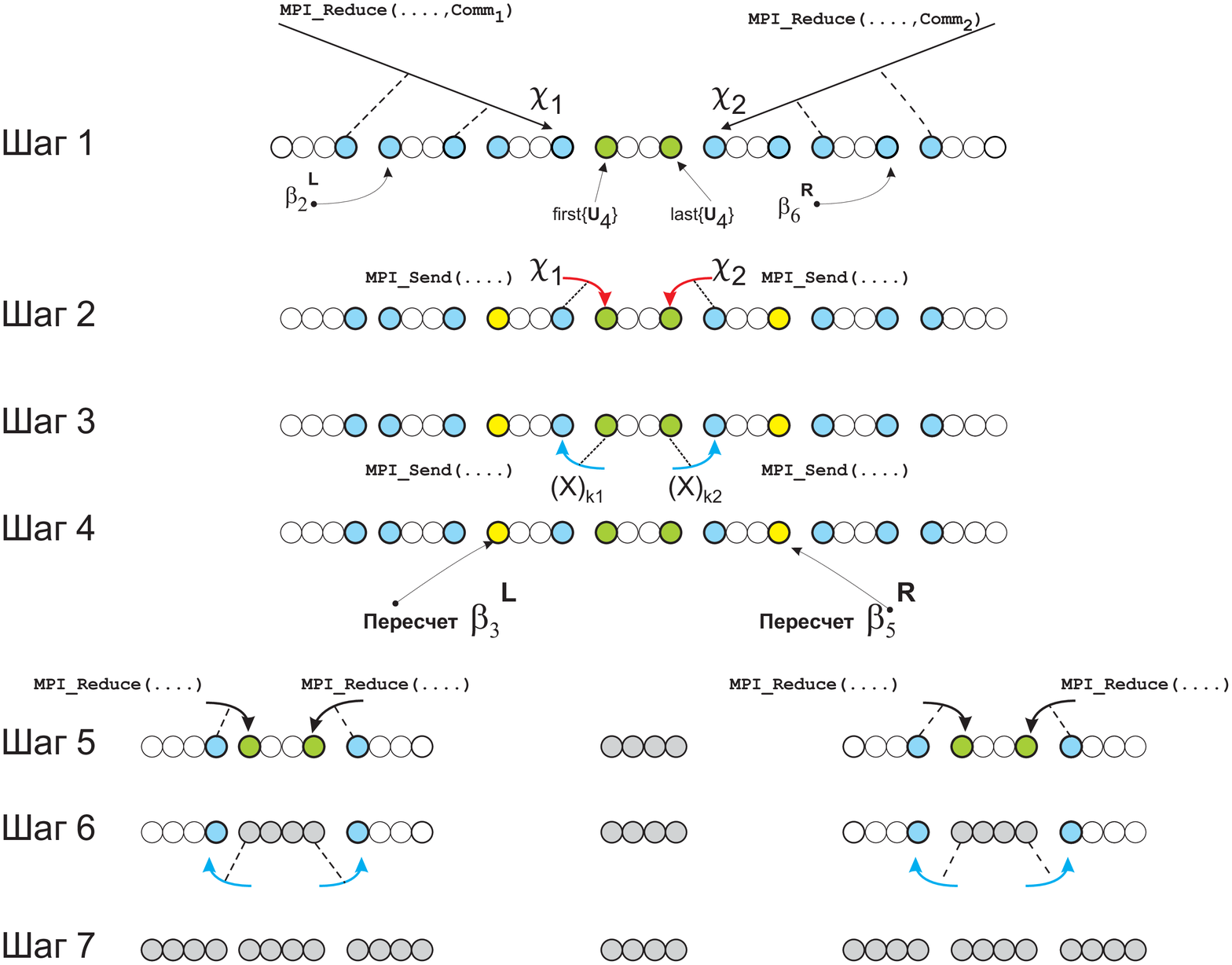}
\caption{Этап разделения СЛАУ в случае $p=7$.}
 \label{pic:scheme}
\end{figure}

Далее, на каждом уровне $s=1,...,\lceil\log_2p\rceil$ процесса
дихотомии из (рис.~\ref{pic:scheme})  вычисляются $2^{s-1}$
элементов из вектора решения по формуле (\ref{theor_3}). В данном
случае на первом уровне дихотомии вычисляются компоненты
\mbox{$\left(\mathbf{X}\right)_{k_1}=first\{\mathbf{U}_4\}$} и
\mbox{$\left(\mathbf{X}\right)_{k_2}=last\{\mathbf{U}_4\}$} .

На шаге $1$ посредством вызова коллективной функций MPI\_Reduce
над коммуникаторами $\mathrm{Comm_1},\mathrm{Comm_2}$ вычисляются
величины

\begin{equation}
\label{step21}
\begin{array}{ll}\displaystyle
\chi_1=\sum_{j=1}^{3}
\beta_{j}^\mathrm{R}\left(\mathbf{Z}_j^\mathrm{R}\right)_{k_1}, &
\displaystyle \chi_2=\sum_{j=5}^{7}
\beta^\mathrm{L}_{j}\left(\mathbf{Z}_j^\mathrm{L}\right)_{k_2}.
\end{array}
\end{equation}

На шаге $2$ процессор $3$ посылает процессору $4$ величину
$\chi_1$, а процессор $5$ -- величину $\chi_2$

Теперь необходимые компоненты, расположенные на $4$ --ом
процессоре могут быть вычислены, как

\begin{equation}
\label{theor33}
\begin{array}{ll} \displaystyle
\left(\mathbf{X}\right)_{k_1}=first\{\mathbf{U}_4\}=\chi_1+\chi_2\frac{\left(\mathbf{G}^{\mathrm{L}}_4\right)_{k_2}}{\left(\mathbf{G}^{\mathrm{R}}_4\right)_{k_2}}+\beta_{4}^{\mathrm{L}},\\\\
\displaystyle
\left(\mathbf{X}\right)_{k_2}=last\{\mathbf{U}_4\}=\chi_1\frac{\left(\mathbf{G}^{\mathrm{R}}_4\right)_{k_1}}{\left(\mathbf{G}^{\mathrm{L}}_4\right)_{k_1}}+\chi_2+\beta_{4}^{\mathrm{R}}.
\end{array}
\label{step23}
\end{equation}

На шаге $3$, для того чтобы исключить найденные компоненты из
системы уравнений (произвести разделение систем), процессор с
номером $4$ посылает процессору $3$ величину
$\delta^{\mathrm{L}}=\left(\chi_{2}\frac{\left(\mathbf{G}^{\mathrm{L}}_4\right)_{k_2}}{\left(\mathbf{G}^{\mathrm{R}}_4\right)_{k_2}}+\beta_{4}^{\mathrm{L}}\right)\left(\mathbf{Z}_{4}^\mathrm{L}\right)_{k_1-1}$ , а процессору $5$ величину
$\delta^{\mathrm{R}}=\left(\chi_{1}\frac{\left(\mathbf{G}^{\mathrm{R}}_4\right)_{k_1}}{\left(\mathbf{G}^{\mathrm{L}}_4\right)_{k_1}}+\beta_{4}^{\mathrm{R}}\right)\left(\mathbf{Z}_{4}^\mathrm{R}\right)_{k_2+1}$ соответственно.

На шаге 4 модифицируется вектор правой части СЛАУ
(\ref{systema}),(\ref{systemb}),(\ref{systemc}); как следствие
(\ref{beta_symn}), на процессорах с номерами 3 и 5 перевычисляются
величины

\begin{equation}
\label{recalc_beta}
\begin{array}{cc}
\hat{\beta}_{3}^\mathrm{R}=\beta_{3}^\mathrm{R}+\delta^{\mathrm{L}},&
\hat{\beta}_{5}^\mathrm{L}=\beta_{5}^\mathrm{L}+\delta^{\mathrm{R}},\\\\
\hat{\beta}_{3}^\mathrm{L}=\beta_{3}^\mathrm{L}+\delta^{\mathrm{L}}\frac{\left(\mathbf{G}^{\mathrm{L}}_{3}\right)_{t_R}}{\left(\mathbf{G}^{\mathrm{R}}_{3}\right)_{t_R}},&
\hat{\beta}_{5}^\mathrm{R}=\beta_{5}^\mathrm{R}+\delta^{\mathrm{R}}\frac{\left(\mathbf{G}^{\mathrm{R}}_{5}\right)_{q_L}}{\left(\mathbf{G}^{\mathrm{L}}_{5}\right)_{q_L}},\\\\

\end{array}
\end{equation}
где $\left(\mathbf{X}\right)_{t_R}\equiv
last\{\mathbf{U}_{3}\}, \quad \left(\mathbf{X}\right)_{q_L}\equiv
first\{\mathbf{U}_{5}\},
$
На следующем уровне дихотомии аналогичным образом разделяются
системы (\ref{systema}),(\ref{systemc}). Для этого в соответствии
с (\ref{theor_3}) вычисляются решения для компонент
$first\{\mathbf{U}_2,\mathbf{U}_6\}$,
$last\{\mathbf{U}_2,\mathbf{U}_6\}$. После завершения процесса
разделения исходная система будет разделена на независимые
подсистемы (\ref{systemb}), для решения которых может быть
применен какой-либо из вариантов метода
прогонки\cite{Ilin_Kuznecov}.

В дальнейшем для определения решения для другой правой части из
(\ref{main_eq3}) потребуется перевычислить только величины
$\beta^{\mathrm{R,L}}$, тогда как векторы
$\mathbf{Z}^{\mathrm{R,L}},\mathbf{G}^{\mathrm{R,L}}$ постоянны
для всех правых частей, так как зависят только от вида матрицы.

\subsection{Вычислительные и коммуникационные затраты процесса дихотомии.}
В \cite{Terekhov:dichotomy} показано, что при реализации операции
All-to-One-Reduce на основе алгоритма \cite{Rabenseifner} оценка
времени процесса дихотомии имеет вид
\begin{equation}
\begin{array}{l}
\ddot{T}^{Dichotomy}_{p}=\sum_{i=1}^{\log_2(p)}\left(2\log_2(2^i)\alpha+\frac{2^i-1}{2^i}l(2\gamma+\beta)\right)=\\\\
=\alpha\left[\log_2(p)+1\right]\log_2(p)+2l\left(\log_{2}(p)-\frac{p-1}{p}\right)\left(\gamma+\beta/2\right),
\end{array}
\label{dich_time2}
\end{equation} где $\alpha$ -- время задержки до начала передачи
данных между процессорами (латентность), $\beta$ -- время передачи
одного числа, $\gamma$ -- время выполнения операции сложения двух
чисел, $l$ -- число серии одновременно решаемых СЛАУ. Сравнивая
полученную оценку с оценкой для алгоритма циклической редукции
\cite{hockney:cyclic}

\begin{equation}
\ddot{T}^{Cyclic\,Reduction}_{p}=2\log_2(p)\left(\alpha+l\beta+l\gamma\right),
\label{cyclic_time2}
\end{equation}

заключаем, что алгоритм дихотомии формально требует большего
времени счета, однако в случае, если решаются несколько серий
задач и величина латентности $\alpha$ незначительна, то оценки
(\ref{dich_time2}),(\ref{cyclic_time2}) практически совпадают.

Высокая производительность алгоритма дихотомии по сравнению с
методом циклической редукции обеспечивается за счет сокращения
времени коммуникаций и времени синхронизации вычислений.
Сокращение времени коммуникации возможно благодаря тому, что
основная коммуникационная операция "+"\  Алгоритма Дихотомии
обладает свойством ассоциативности. Архитектура современных
суперкомпьютеров такова, что время парных взаимодействий для
различных процессоров может существенно различаться
\cite{Dongarra:Reduce,Voevodin_MPI,MPI_REDUCTION}. Ассоциативность
вычислений разрешает на уровне коммуникационной библиотеки или
языка программирования определить порядок взаимодействий ПЭ таким
образом, чтобы учитывая архитектуру суперкомпьютера,
минимизировать время обменов данными. Организация обменов через
вызовы неблокирующей функции "All-to-One-Reduce(+)" позволяет
сократить время синхронизации процессорных элементов.
Действительно, если в рамках одной группы процессоров
\footnote{Коммуникатор -- в терминологии MPI.} существуют два
свободных ПЭ с подготовленными данными, то коллективная операция
"(+)"\ над этой группой процессоров может начинать выполняться,
даже если предыдущий вызов "All-to-one-Reduce"\ на всех
процессорах не завершен. Например, процессоры с номерами $1$ и $2$
(см. рис.~\ref{pic:scheme}) на шаге 1 могут сложить/передать свои
данные и начать выполнение шага $5$ без синхронизации вычислений с
другими процессорами.

\section{Численные эксперименты}
Довольно часто эффективные с теоретической точки зрения
параллельные алгоритмы при их реализации на суперкомпьютерах могут
не обеспечивать ожидаемого сокращения времени счета. В первую
очередь это связано с тем, что при анализе эффективности
конкретного алгоритма достаточно сложно объективно учесть все
особенности вычислительных систем (время доступа к памяти,
латентность и пропускную способность сети и т.д.). Поэтому
численные эксперименты с модельными постановками задач являются
важным этапом исследования параллельных алгоритмов.

\subsection{Решение уравнения Акустики.}
В цилиндрической системе координат $(r,z)$ в полупространстве
$z\geq 0$ рассмотрим проблему моделирования распространения
акустических волн от точечного источника

\begin{equation}
\label{acoustic_problem}
\begin{array}{llr}
\displaystyle \frac{1}{\rho^2({\bf x} )}\frac{\partial^2
u}{\partial t^2}({\bf x},t)=\nabla \left[V_s({\bf x})\, \nabla
u({\bf x},t)\right]+\delta({\bf x-x_0})f(t),& t>0,\quad {\bf
x}=(r,z).
\end{array}
\label{wave-eq}
\end{equation}

Будем полагать, что проблема (\ref{acoustic_problem}) решается при
нулевых начальных данных
\begin{equation}
\begin{array}{llr}
\displaystyle \left. u\right|_{t=0}=\left.\frac{\partial
u}{\partial t}\right|_{t=0}=0.
\end{array}
\label{initcond}
\end{equation}
Будем считать, что при $z=0$ поверхность является свободной, а по
координате $r$ и $z$ введены вспомогательные границы
\begin{equation}
\begin{array}{llr}
\displaystyle \left.\frac{\partial u}{\partial
z}\right|_{z=0,l_2}=\left.u\right|_{r=l_1}=0.
\end{array}
\label{boundary_cond}
\end{equation}
Границы $r=l_1$ и $z=l_2$ выбираются таким образом, чтобы для
рассчитываемого момента времени не возникало волн, отраженных от
них. Дополнительно потребуем выполнения условия
\cite{Samarski_andreev}
\begin{equation} \left.\frac{\partial u}{\partial
r}\right|_{r=0}=0. \label{symm_cond}
\end{equation}

Будем искать решения задачи
(\ref{acoustic_problem})-(\ref{symm_cond}) в виде ряда Фурье по
функциям Лагерра \cite{Mikh:Laguerre}

\begin{equation}
Q_m(\mathbf{x})=\int_{0}^{\infty}u(\mathbf{x},t)(ht)^{-\frac{\alpha}{2}}l_{m}^{\alpha}(ht)dt
\label{series_lag}, \quad
u(\mathbf{x},t)=(ht)^{\frac{\alpha}{2}}\sum_{m=0}^{\infty}Q_m({\bf
x})l^{\alpha}_m(ht) .
\end{equation}

где  $l^{\alpha}_m(ht)$ ортонормированные функции Лагерра
\cite{Suetin}, которые выражаются через классические полиномы
Лагерра следующим образом

$$
l^{\alpha}_m(ht)=\sqrt{\frac{hm!}{(m+\alpha)!}}(ht)^{\frac{\alpha}{2}}e^{-\frac{ht}{2}}L^{\alpha}_m(ht).
$$
Здесь $m$ -- степень полинома Лагерра, $\alpha \ge 2$ -- целая
константа, $h$ -- параметр преобразования. Применяя преобразование
(\ref{series_lag}) к (\ref{acoustic_problem})--(\ref{symm_cond}),
получаем краевую задачу для определения $m$-го коэффициента
разложения

\begin{equation}
\left\{
\begin{array}{l}
\label{laguerre_h} \displaystyle \nabla \left[V_s({\bf x})\,
\nabla Q_m({\bf x})\right]-\frac{h^2}{4\rho^2({\bf x})}Q_m({\bf
x})=-\delta({\bf
x-x_0})f_m+\frac{h^2}{\rho^2({\bf x})}\sqrt{\frac{m!}{(m+\alpha)!}}\sum_{k=0}^{m-1}(m-k)\sqrt{\frac{(k+\alpha)!}{k!}}Q_k({\bf x}),\\\\
\displaystyle\left. \frac{\partial Q_m}{\partial
r}\right|_{r=0}=\left.\frac{\partial Q_m}{\partial
z}\right|_{z=0,l_2}=\left.Q_m\right|_{r=l_1}=0,
\end{array}\right.
\end{equation}
где
$f_m=\int_0^{\infty}f(t)(ht)^{-\frac{\alpha}{2}}l^{\alpha}_{m}(ht)dt$.

Данный метод можно рассматривать как аналог
спектрально-разностного метода на основе Фурье-преобразования
\cite{FEM1}, однако, здесь роль "частоты"\ выполняет параметр $m$,
определяющий степень полиномов. Другое отличие преобразования
Лагерра от Фурье состоит в том, что параметр разделения гармоник
присутствует только в правой части (\ref{laguerre_h}). Таким
образом, чтобы рассчитать динамику волнового поля необходимо
решить серию уравнений вида (\ref{main_problem}). Такая постановка
позволяет использовать Алгоритм Дихотомии в рамках
предобуславливающей процедуры, рассмотренной выше.

\subsection{Обращение предобуславливающего оператора.}
В рамках схемы (\ref{main-scheme}) для решения задачи
(\ref{laguerre_h}) выберем в качестве предобуславливающего
оператор

\begin{equation}
B\equiv \Lambda_r+\Lambda_z-d, \label{precond}
\end{equation}

где $\Lambda_r,\Lambda_z$ определены в (\ref{scheme_op_r_z}), а
коэффициенты полагаются следующими\footnote{Введем обозначение
$\tilde{f}=\frac{1}{2}\left(\min_{x \in G }f({\bf x})+\max_{x \in
G }f({\bf x})\right).$}

\begin{equation}
\label{diff_approx2}
\begin{array}{ccc}
 \displaystyle
a_1(i,k)=\bar{r}_i\tilde{V},& a_2(i,k)=r_i\tilde{V},&
\displaystyle d(i,k)=r_i\frac{h^2}{4}\widetilde{\rho^{-1}}. \\
\end{array}
\end{equation}

Для решения задачи $By=f$ воспользуемся методом разделения
переменных\cite{Samarski_Nikolaev,FFT_Poisson}. Поскольку
\mbox{$\left.\frac{\partial y}{\partial z}\right|_{z=0,l_2}=0$},
будем искать решение в виде ряда по собственным функциям
разностного оператора $\Lambda_z $:

\begin{equation}
\begin{array}{ll}
\displaystyle
y_{i,k}=\sqrt{\frac{2}{N_2}}\left(\frac{1}{2}\tilde{y}_i(1)+\sum_{l=2}^{N_2}\tilde{y}_i(l)\cos\left(\frac{\pi
(k-1/2)(l-1)}{N_2}\right)\right), & 1\leq i\leq N_1-1,\quad  1\leq
k \leq N_2,\\\\
\displaystyle
\tilde{\phi}_{i}(l)=\sqrt{\frac{2}{N_2}}\sum_{k=1}^{N_2}f_{i,k}\cos\left(\frac{\pi(k-1/2)(l-1)}{N_2}\right),&
1\leq i\leq N_1-1,\quad 1\leq l \leq N_2.
\end{array}
 \label{fft_phi}
\end{equation}
где $\tilde{y}_{i}(l)$ определяются из решения трехдиагональной
системы уравнений

\begin{equation}
\begin{array}{cr}
\displaystyle \label{eq_trid}
\left(\Lambda_r+4\frac{a_2}{h_z^2}\sin^2\frac{\pi
(l-1)}{2N_2}-d\right)\tilde{y}(l)=\tilde{\phi}(l),& l=1,...,N_2.
\end{array}
\end{equation}

Суммы (\ref{fft_phi}) будем вычислять, используя последовательный
алгоритм быстрого дискретного преобразования Фурье
\cite{Samarski_Nikolaev,FFT_TURK}, а для определения решений из
серии уравнений (\ref{eq_trid}) -- использовать Алгоритм Дихотомии
для решения серии трехдиагональных систем уравнений.

Остановимся на некоторых аспектах реализации Алгоритма Дихотомии в
рамках метода разделения переменных. Одно из ограничений на
область применимости Алгоритма Дихотомии состоит в том, что все
СЛАУ из серии должны включать одну и ту же фиксированную матрицу
(\ref{main_eq3}). Это позволяет пренебрегать затратами на
подготовительные вычисления. Однако для (\ref{eq_trid}) условие
постоянства матрицы не выполняется, поэтому для этого случая
Алгоритм Дихотомии не применим. С другой стороны, необходимость
многократного обращения оператора позволяет рассматривать
совокупность уравнений (\ref{eq_trid}) как серию задач вида
(\ref{main_eq3}) для фиксированного $l$ и, следовательно,
использовать Алгоритм Дихотомии.

\subsection{Оценка производительности.}
\begin{figure}[!h]
\begin{center}
\includegraphics[width=\textwidth]{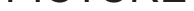}
\end{center}
\caption{Модель среды "Геологический сброс"\ .} \label{pic:Sbros}
\end{figure}

В качестве тестовой выберем модель среды типа "Геологический
сброс"\ , где скорость распространения волн в среде
$V_s(\mathbf{x})$ задана в соответствии с рис.~\ref{pic:Sbros}.a,
при этом плотность среды будем полагать постоянной
$\rho(\mathbf{x}) \equiv 1$.

Источник волнового поля расположим в начале координат, а временную
зависимость сигнала зададим в виде
\begin{equation}
\label{functiont}
 f(t)=\exp\left[-\frac{(2\pi
f_0(t-t_0))^2}{\gamma^2}\right]\sin(2\pi f_0(t-t_0)),
\end{equation}
с параметрами $f_0=30\mathrm{Hz},\;t_0=0.2s,\;\gamma=4$.
Аппроксимация уравнения (\ref{laguerre_h}) проводилась  на
равномерной сетке $\bar{\omega}$ c числом узлов
\mbox{$N_1=N_2=2^k$}, \mbox{$k=\{11,12,13,14,15\}$}. Число
слагаемых в ряде (\ref{series_lag}) было ограничено $n=2000$, а
параметры разложения полагались следующими $\alpha=5,\;h=300$.

\begin{figure}[!h]
\begin{center}
\includegraphics[width=0.45\textwidth]{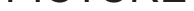}\hfill
\includegraphics[width=0.5\textwidth]{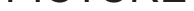}
\parbox[t]{0.45\textwidth}{ \caption {Декомпозиция расчетной области}\label{pic:decomposition} }
\hfill
\parbox[t]{0.45\textwidth}{\caption{Зависимость времени счета от числа процессоров для сетки $N_1=N_2=2^k$,где $k=$  a)$\,12$; b)$\,13$; c)$\,14$; d)$\,15$. }\label{main_pic23}}
\end{center}

\end{figure}

Для оценки производительности предлагаемого алгоритма были
реализованы неявный трехслойный Чебышевский и CG-методы, а для
обращения предобуславливающего оператора -- метод разделения
переменных. Начальное приближение в схеме (\ref{main-scheme}) для
всех $Q^{0}_n,\; n=1,...,2000$ задавалось нулевым во всей
расчетной области, а итерационный процесс завершался по достижению
условия
$$\displaystyle \frac{\parallel AQ^{k}_{n}-f\parallel}{\parallel f
\parallel} \leq 10^{-10},$$ где $k$--номер итерации.
С учетом выбранной декомпозиции расчетной области
рис.~\ref{pic:decomposition},  решение трехдиагональных систем
уравнений (\ref{eq_trid}) выполнялось в направлении $r$, а
Фурье-преобразование (\ref{fft_phi}) проводилось в направлении
$z$.

Все численные процедуры были реализованы на языке Fortran-90 с
использованием MPI-технологии. Для выполнения быстрого
Фурье-преобразования использовалась библиотека FFTW\cite{FFTW}.
Вычисления проводились на суперкомпьютере "МВС-100k"\
Межведомственного Суперкопьютерного Центра РАН, построенного на
основе четырехядерных процессоров Intel Xeon, работающих на
частоте 3 ГГц и соединенных коммуникационной средой Infiniband.

На рис.~\ref{main_pic23}. и табл.~1 приведены зависимость времени
счета ($\mathbf{T}$) и величины ускорения  ($\mathbf{S}$) от числа
процессоров для одной итерации при решении задачи
(\ref{laguerre_h}) методом Чебышева и методом сопряженных
градиентов. При реализации этих алгоритмов на суперкомпьютере было
достигнуто почти линейное ускорение для сеток различной
подробности в широком диапазоне числа процессоров. Достигнутый
уровень ускорения обеспечивается за счет применения Алгоритма
Дихотомии, так как основные коммуникационные затраты приходятся на
решение трехдиагональных СЛАУ при обращении предобуславливающего
оператора. Таким образом, требования обеспечения высокой
массштабируемости алгоритма и возможности проведения расчетов с
использованием тысяч процессоров можно считать выполненными.

\begin{table}[!h]
\center \small
\begin{tabular}{lcccccccccc}
  \hline
  size & \multicolumn{2}{c}{2048x2048} &\multicolumn{2}{c}{4096x4096}& \multicolumn{2}{c}{8192x8192}&\multicolumn{2}{c}{16384x16384}&\multicolumn{2}{c}{32768x32768} \\ \hline
   NP & $ \mathrm{T}$&  $ \mathrm{S}$  &  $ \mathrm{T}$&  $ \mathrm{S}$  & $ \mathrm{T}$&$ \mathrm{S}$&$ \mathrm{T}$&$ \mathrm{S}$&$ \mathrm{T}$&$ \mathrm{S}$\\
  \hline
  64 &1.4e-02 &-&   8.0e-02&     - & 3.5e-01   & -&  1.7      &-&-&-\\
  128 & 7.3e-03&122& 3.7e-02  &138&   1.8e-01   &124&   7.2e-01    &151&3.9&-\\
  256 & 6.3e-03&142& 1.8e-02  &284&   9.0e-02   &254&   4.3e-01    &253&2.15&172\\
  512 & -&-& 1.1e-02  & 465 & 5.0e-02   &448&   2.0e-01     &544&1.01&463 \\
  1024 &-&-& 1.0e-02  & 512 & 2.7e-02   &829&       1.0e-01&1088&5.4e-01&924\\
  2048 &- &-&-        &-    & 2.0e-02   &1120&7.0e-02&1554&3.2e-01&1560\\
  \hline
\end{tabular}
\caption{Зависимость времени счета одной итерации ($\mathrm{T}$) и
коэффициента ускорения ($S$) от числа процессоров для CG-метода.}
\end{table}

Выше обращалось внимание на то, что из-за наличия коллективных
взаимодействия типа "All-to-All-Reduce"\ ,  реализация метода
сопряженных градиентов на суперкомпьютере может быть менее
эффективна, чем использование Чебышевского набора параметров.
Однако вычислительные эксперименты не выявили существенных
различий в зависимости величины ускорения от числа процессоров для
этих методов(рис.~\ref{main_pic23}) . Действительно, сопоставляя
оценки времени коммуникационных взаимодействий для CG-Метода и
Алгоритма Дихотомии
$$
\begin{array}{l}
T^{CG,\,all-reduce}_{p}=2\log_2(p)\alpha+\frac{p-1}{p}\left(\gamma+2\beta\right),\\\\
T^{Dichotomy}_{p}=\alpha\left[\log_2(p)+1\right]\log_2(p)+l\left(\log_{2}(p)-\frac{p-1}{p}\right)\left(\gamma+2\beta\right),
\end{array}
$$
заключаем, что для вычислительных систем с низкой величиной
латентности $\alpha$ при условии \mbox{$l\gg 1$}, время
межпроцессорных взаимодействий CG-метода на фоне коммуникационных
взаимодействий Алгоритма Дихотомии незначительно.

\begin{figure}[!h]
\begin{center}
\includegraphics[width=0.5\textwidth]{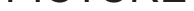}
\end{center}
\caption{Модель среды Marmousi.} \label{pic:marmous}
\end{figure}

С целью оценки точности получаемого решения, для модели
"Геологический Сброс"\ на рис.~\ref{pic:Sbros}.a  приведено
волновое поле момент времени $t=2.5$. На рис.~\ref{pic:Sbros}.b.
для сеток различной подробности приведена зависимость амплитуды
волнового поля вдоль прямой $Slice(z)=u(1.2km,z)$. Видно, что
точность получаемого решения значительно зависит от числа узлов
сетки, приходящееся на характерную длину волны. Принимая, решение,
полученное с использованием самой подробной сетки
$h_r=h_z=1/280\lambda,\; \lambda\simeq33m$ в качестве точного,
получаем следующие оценки скорости сходимости численной схемы:

\begin{table}[!h]
\center \small
\begin{tabular}{c|c|c|c|c}
\hline $h_{r,z}$& $N_\lambda$&$\frac{\parallel
\tilde{u}-u\parallel}{
\parallel u\parallel}$&$\frac{\parallel
u-\tilde{u}_{1/2}\parallel}{
\parallel u-\tilde{u}_{1/4}\parallel}$&$M_{\Delta}$
  \\
  \hline
  % after \\: \hline or \cline{col1-col2} \cline{col3-col4} ...
  $1/35\lambda$ &  $4096\times4096$&$0.45$& -&90840\\
  \hline
  $1/70\lambda$ &  $8192\times8192$ &$0.11$&4.09&90590\\
  \hline
  $1/140\lambda$ & $16384\times16384$&$0.022$&5&90680 \\
  \hline
\end{tabular}
\caption{Зависимость точности решения и количества обращений
предобуславливающего оператора от числа узлов сетки.}
\label{tabl:sbros}
\end{table}

Видно, что при уменьшении пространственного шага в два раза
величина погрешности в сеточной норме $L_2$ уменьшается в $4$ и
$5$ раз соответственно. Таким образом, предлагаемый алгоритм
обеспечивает второй порядок точности по пространству.
Дополнительно отметим, что число обращений предобуславливающего
оператора $M_\Delta$ при решении задачи (\ref{wave-eq})
практически не зависит от числа узлов сетки. Таким образом,
предлагаемая предобуславливающая процедура эффективна в случае,
когда требуется проводить расчеты с высоким пространственный
разрешением при умеренной контрастности среды.

\begin{figure}[!h]
\begin{center}
\includegraphics[width=0.5\textwidth]{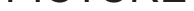}\hfill
\includegraphics[width=0.5\textwidth]{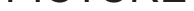}
\parbox[t]{0.45\textwidth}{ \center {a) Приемник $\mathbf{x}_1=(30\lambda_1,0)$ }}\hfill
\parbox[t]{0.45\textwidth}{ \center {б) Приемник $\mathbf{x}_2=(70\lambda_1,0)$} }
\end{center}
\caption{Сейсмограммы для среды Marmousi.}
\label{pic:trace_marmousi}
\end{figure}

Дополнительно для реальных пространственно-временных масштабов
было произведено моделирование распространения акустических волн
для модели среды "Marmousi"\ (рис.~\ref{pic:marmous}a,б)
\cite{Marmousi}. На рис.~\ref{pic:trace_marmousi} приведена
зависимость амплитуды поля от времени для двух приемников с
координатами
$\mathbf{x}_1=(30\lambda_1,0),\;\mathbf{x}_2=(70\lambda_1,0)$,
источник вида (\ref{functiont}) так же был помещен в начале
координат. Расчеты проводились на сетках с числом узлов $N_r\times
N_z=\{4096\times1321,8192\times2642,16384\times5284,32768\times10568\}$,
что соответствовало пространственным шагам
$h_r=h_z=\{1/25\lambda_1,1/50\lambda_1,1/100\lambda_1,1/200\lambda_1\}$,
$\lambda_1\simeq50m$. Число обращений предобуславливающего
оператора для определения всех гармоник $Q_n,\;n=1,...,2000$ было
порядка $90000$ и практически не зависело от шага сетки.

Так согласно рис.~\ref{pic:trace_marmousi}.a и
рис.~\ref{pic:trace_marmousi}.b разрешающая способность сетки
$1/25\lambda_1,\,1/50\lambda_1$ не обеспечивает приемлемой
точности ни для первого, ни для второго источника. При разрешающей
способности сетки в $1/100\lambda_1$ точность расчетов для
приемника $\mathbf{x}_1$ становится достаточно высокой, однако,
для приемника $\mathbf{x}_2$, который фиксирует волну в более
поздний момент времени, точность расчетов несколько меньше. Этот
эффект обусловлен так называемой "фазовой ошибкой"\ характерной
для численно-аналитических методов\cite{FEM1}. Таким образом, при
моделировании волновых процессов для длительных временных периодов
необходимо использовать сетки с высокой разрешающей способностью,
чтобы погрешность расчетов для заключительных моментов времени
оставалась допустимой.

\section{Заключение.}
В работе предложен параллельный алгоритм для решения эллиптических
уравнений второго порядка c неразделяемыми переменными. Разностные
уравнения, полученные в результате конечно-объемной аппроксимации,
решаются неявным полуитерационным методом Чебышева или CG-методом.
Выбор оператора Лапласа в качестве предобуславливающего позволяет
обеспечить высокую скорость сходимости итерационного процесса для
сред с умеренной контрастностью, при этом число итераций для
достижения заданной точности практически не зависит от числа узлов
сетки.

Почти линейная зависимость коэффициента ускорения и высокая
масштабируемость параллельного алгоритма обеспечиваются за счет
использования Алгоритма Дихотомии в рамках метода разделения
переменных для обращения предобуславливающего оператора. Таким
образом, высокая производительности Алгоритма Дихотомии и простота
его реализации позволяет эффективно распараллеливать экономичные
численные процедуры, требующие многократного решения
трехдиагональных систем уравнений.

В рамках вычислительных экспериментов установлено, что при
моделировании распространения звуковых волн в неоднородной среде
для реально пространственно-временных масштабов необходимо
использовать сетки с высокой разрешающей способностью. Таким
образом, предложенный алгоритм является востребованным, так как
позволяет с высокой эффективностью и в полной мере использовать
вычислительные мощности современных суперкомпьютеров при решении
больших задач.

\section{Благодарности.}
Автор благодарит Фатьянова А.Г. и Малышкина В.Э. за многочисленные
обсуждения результатов данной работы.

 \end{document}